
\magnification=\magstep1
\input amstex
\documentstyle{amsppt}
\pageheight{8.5truein}
\NoBlackBoxes\leftheadtext{Propagation of Localization and Optimal Entropy}
\rightheadtext{E. Carlen and A. Soffer}

\topmatter
\title Propagation of Localization Optimal Entropy production and
convergence rates for the central limit Theorem
\endtitle
\author E. Carlen and A. Soffer\endauthor
\address Mathematics Department,
 Rutgers University, New Brunswick, NJ 08903 \endaddress
\address Mathematics Department,
 Rutgers University, New Brunswick, NJ 08903
\endaddress
\email soffer\@math.rutgers.edu
\endemail
\abstract We prove for the rescaled convolution map $f\to
f\circledast f$ propagation of polynomial, exponential and gaussian
localization.  The gaussian localization is then used to prove an
optimal bound on the rate of entropy production by this map.  As an
application we prove the convergence of the CLT to be at the optimal
rate $1/\sqrt{n}$ in the entropy (and $L^1$) sense, for
distributions with finite 4th moment.
\endabstract
\endtopmatter

\document
\subhead
Section 1. - Introduction, Notation, Preliminaries\endsubhead
\medskip

The Central limit Theorem (CLT) naturally leads to the analysis of the
(nonlinear) rescaled convolution map, of a probability density with itself.  Related maps
appear in the study of Boltzmann type equations.  A major issue is the
convergence and rate in various norms for CLT.  In this work, we will
study the convergence in the strong norm $L^1$, and the stronger sense
of convergence in relative entropy.

To find rate, we use monotonicity or
entropy production estimates for the convolution map convergence in
this sense was first established by  Barron [Bar].  The corresponding
result for the Boltzman equation was established by Carlen, Carvalho
and Wennberg [CCW].  Such estimates
have also allowed, via the method of [CS] to prove the CLT for
dependent variables, in a nonperturbative way.

Our main tool is  an optimal entropy production rate for the
convolution map;  such estimate depends critically on {\bf
propagation of localization};   to successfully apply then entropy
production bound, one needs to show that the localization at
infinity is not spoiled under iteration of  the convolution map.  We
prove in sections 2 and 3 that polynomial exponential and,  most
importantly, gaussian localization are uniformly propagated  the
convolution map.  These results are then used to derive the optimal
entropy production bounds in the gaussian case, and as application
gives the optimal $1/\sqrt{n}$ convergence of the CLT in the
entropy, and $L^1$ norms, for gaussians (or better) localization, as
well as the case of bounded moments to order 4.

Propagations of localization are important for other applications.  For
example, gaussian propagation of localization for the Boltzmann kernel
would have major implications to asymptotic stability
and more. [CC1,2, CGT, CELMR, Des, De94, GTN]

We conclude with some mention of possible applications.  Our proof of
the propagation of localization in the polynomial and exponential
cases is rather direct.  In the polynomial case it follows from moment
estimates and in the exponential case by direct estimates on the
generating function.

The Gaussian case is however much more difficult.  It is based on a
kind of asymptotic log concavity in the CLT, combined with a theorem of Brascamp
and Lieb, and other analytic arguments. The estimates of entropy
production uses linear approximation theory of the map, combined
with gaussian propagation of localization to arrive at the leading
entropy growth term. The propagation of gaussian localization, which
is crucial for getting the optimal convergence rate for the CLT, is
based on upper AND lower bounds on the distribution $\rho.$ Hence,
if the distribution has a thin tail, it results in delocalization of
the entropy, which breaks the needed estimates. This problem is
usually overcome by assuming, on top of the localization, a spectral
gap assumption [BaBN,Bart,Jon,Vil  ].

We use a new construction to overcome this problem, thus avoiding
the assumption of spectral gap, and extending the optimal
convergence rates to arbitrarily gaussian localized distribution,
with finite Fisher information.

As we shall show, if a density $\rho$  has most of its mass localized in the sense of having sufficiently many
moments bounded, and if we are given a bound on the Fisher information of $\rho$, then the tails of
$\rho$ do not contribute significantly to the the entropy of $\rho$, not to the entropy production by
rescaled contribution of $\rho$. Without the bound on the Fisher information, this would not be the
case at all. But since abounds on Fisher information are rescaled by iterated convolution, this opens
the way to the following strategy for dealing with possibly thin tails:
{\it We approximate $\rho_n$ by a new
distribution, $\tilde \rho_n,$ which is obtained by stitching a
gaussian tail to $\rho$, for $|x|\ge c \sqrt n$, and renormalizing
the mean and variance.}
Then, we show that the monotonicity
estimates are optimal for the stitched distribution, and the
difference to $\rho$ is exponentially small. The effect of the small
errors is absorbed by the monotonicity (entropy production) bounds,
similar to the way perturbations of the convolution map were treated
in our paper [CS].

Our notation and preliminaries follow closely the paper [CS].  Here we
briefly recall the main ingredients of entropy/information
bounds. [CS, Dem, Lie78, Lie89, Bar]

Let $X$ be an $\Bbb R^m$ valued random variable on some probability
space.  Let $\mu$ denote the law of $X$.  If $d\mu(x)= \rho (x) dx$,
we say that $X$ has density $\rho(x)$.  $m(x)$ stands for the mean
of $x$, and $m_j(x)$ for the $j$ the moment of $X$.  The variance is
then
$$
[\sigma(x)]^2 = E(|X-m(x)|^2)
$$
and $X$ has variance 1 if $[\sigma(x)]$ is the identity matrix.  Let $g_t$
denote the centered Gaussian density with variance $t$:
$$
g_t(x) = (2\pi t)^{-m/2}e^{-x^2/2t}.
$$
$$
g\equiv g_1.
$$

The entropy of $\rho$ is
$$
S(\rho) = - \int\rho  \, \ln \rho dx
$$
and the relative entropy of $\rho$ is
$$
D(\rho) = \int\frac{\rho(x)}{g(x)}\left( \ln\frac{\rho(x)}{g(x)}\right)g(x) dx.
$$
By Jensen's inequality $D(\rho)\geq 0$ with equality just when $\rho =
g$.  Clearly, if $\rho$ has mean zero and unit variance
$$
-\infty\leq S(\rho) \leq S(g)
$$
and the upper bound is saturated only when $\rho=g$.

Moreover, for $\rho$ with mean zero and unit variance, which we will
refer to as $\rho$ being {\bf normalized},
$$
D(\rho) = S(g)-S(\rho).
$$
For centered density $\rho$ with $\sigma^2(\rho) = T_r[\sigma(X)^2]$
($m$- the dimension) and $\sqrt{\rho} \in H^1(\Bbb R^m)$, the
Sobolev space, we define the Fisher information
$$
I(\rho) = 4\int_{\Bbb R^m} |\nabla\sqrt{\rho(x)}|^2 dx
$$
and the relative Fisher information, $J(\rho)$ as
$$
J(\rho) = 4\int|(\nabla + \frac{x}{2})\sqrt{\rho(x)}|^2 dx
$$
Clearly, $J(\rho)\geq0 \text{ with } J(\rho) =0 \iff \rho=g$.

Also, note that, when $J (\rho) < \infty$,
$$
J(\rho) = \int_{\Bbb R^m} |\nabla \ln \rho (x) - \nabla \ln  g(x)|^2
\rho(x) dx.
$$
The origin of the convolution map is the following: Suppose $X_1,
X_2$ are two independent random variables with densities $\rho_1,
\rho_2$.  For $0<\lambda < 1$ denote the density of $\lambda X_1 +
(1-\lambda^2)^{1/2} X_2$ by $\rho_1\underset{1/\lambda}\to{*}
\rho_2$.  One computes
$$
\rho_1\underset{1/\lambda}\to{*} \rho_2(u) = \int_{\Bbb R^m}
\rho_1(\lambda u - (1-\lambda ^2)^{1/2}v)\rho_2((1-\lambda^2)^{1/2}
u + \lambda v) d v.
$$
Let $f$ be a bounded measurable
 function on $\Bbb R^m$. Define the
operator $P_t, t> 0$
$$
P_tf(x) = E f(e^{-t} x + (1-e^{-2t})^{1/2} G)
$$
Then $P_t$ is a contraction semigroup on each $L^p(\Bbb R^m, g(x)
dx)\, 1\leq p \leq \infty$. $P^*_t$ denotes the adjoint in $L^1(\Bbb
R^n, dx)$. In particular, if $X$ is a random variable with density
$\rho$
$$
P_t^*\rho(x)\text{ is the density of }e^{-t} X+(1-e^{-2t})^{1/2} G.
$$
We have the following relation between entropy and information, which
is contained in [CS].

\proclaim{Lemma}  Suppose $\rho$ is a centered density with
$\sigma^2(\rho)$.  Then $t\to S(P^*_t\rho)$ is continuous and
monotone increasing on $[0,\infty)$ with
$$
\lim_{t\to \infty} S(P^*_t\rho) = S(g).
$$
Furthermore, when $S(\rho) > -\infty, t\to S(P^*_t\rho)$ is
continuously differentiable on $(0,\infty)$ and
$$
S(P^*_t\rho) = S(\rho) + \int^t_0 J(P^*_s\rho) ds
$$
and
$$
D(\rho) = \int^\infty_0 J(P^*_t\rho) dt.
$$
\endproclaim

We will also use the inequality
$$
D(x) \leq \frac{1}{2} J(x)
$$
due to Stam [Sta] which is equivalent to Gross's logarithmic Sobolev
inequality [Gro], [Ca].

The proof follows from
$$
D(x) =\int^\infty_0 J(e^{-t} X + (1-e^{-2t})^{1/2} G) dt \leq \int^\infty_0
e^{-2t} J(X) dt.
$$
using the Blackman-Stam inequality:
$$
J(e^{-t} X + (1-e^{-2t})^{1/2} Y) \leq e^{-t}J(X) + (1-e^{-2t})J(Y).
$$
We also have the Kullback-Liebler inequality
$$
\|\rho - g\|^2_{L^1(\Bbb R^m, dx)}\leq 2 D(\rho).
$$
The main inequality we prove for entropy production is that under
favorable assumption on both smoothness and gaussian localization of $\rho$,
$$
S(\rho \underset{\sqrt{2}}\to{*}\rho) - S (\rho) \geq CD(\rho).
$$
Our previous work only gave a lower bound of the form $\Phi_\rho
(J(\rho))$, [CS]. The application of this inequality requires that
localization and smoothness is maintained under repeated iteration.
So, for this we prove that gaussian (polynomial and exponential)
localization is uniform in $n$ for
$$
\rho_n = \rho \underset{\sqrt{2}}\to{*} \rho \cdots
\underset{\sqrt{2}}\to{*} \rho,  \qquad  n\text{ times}.
$$

We now state the {\bf main theorem} with convergence rate:

\proclaim{Theorem  (Optimal Entropy convergence)}  Let $\rho$ be a
regular, normalized, variance 1 and with bounded 4th moment
distribution:
$$
I(\rho)<\infty,
$$

$$
\| \rho |x|^4\|_{1}<c<\infty.
$$

 Then
$$
|D(\rho_N)|\leq c/N ,\tag 5.2
$$
and $N:=2^n$.
In particular, the CLT holds in the Entropy (and $L^1$) sense with
the optimal convergence rate $1/\sqrt{n}$.
\endproclaim

 \subhead
Section 2. Propagation of Localization I - Polynomial and
Exponential\endsubhead

Let $\rho$ be normalized distribution, localized exponentially:
$$
\rho_\alpha (x) \equiv \lambda(\alpha)^{-1} e^{\alpha x}\rho(x)\tag
2.1
$$
with $\rho_\alpha(x)$ bounded and in $L^1$: here $\lambda(\alpha)$ is
the normalization constant so that
$$
\int\rho_\alpha(x) dx = 1.\tag 2.2
$$
Therefore
$$
\lambda(\alpha) = \int e^{\alpha x} \rho(x)\tag 2.3
$$
\proclaim{Theorem 2.1 (Exponential Localization)} Let $\rho$ be a distribution in $L^1$ and such
that $\lambda(\alpha)< \infty$ for $\alpha \leq A, A>0$.

Then
$$
L_{\rho_n}(\alpha) \equiv \int e^{\alpha x}\rho_n(x) dx = \int
e^{\alpha x} \sqrt{n}\rho*\cdots* \rho (\sqrt{n} x) dx \leq 2e^{\alpha^2/2}
\tag 2.4
$$
for all $\alpha < A$.
\endproclaim
\remark{Remark}  The above identity, of equation 2.4, is due to Cramer [Cr].\endremark
\demo{Proof}
First we compute the convolution
$$
\aligned
\rho_\alpha *\rho_\alpha&= \lambda(\alpha)^{-2} \int
e^{\alpha(x-y)}\rho(x-y) e^{\alpha y} \rho(y) dy\\
 &= \lambda(\alpha)^{-2}\int e^{\alpha x}\rho(x-y) \rho(y) dy = (\rho
*\rho)_\alpha.\endaligned\tag 2.5
$$
\enddemo
Therefore since
$$
\rho_n = \sqrt{n}\rho * \rho \cdots* \rho(\sqrt{n} x) , \qquad n\text{ times}
$$
we have
$$
\aligned
L_{\rho_n}(\alpha) &= \int e^{\alpha x}\sqrt{n}\rho *\cdots \rho (\sqrt{n} x)
dx\\
&=\int e^{\frac{\alpha}{\sqrt{n}}y}\rho * \cdots\rho(y) dy\\
&=L^n_\rho(\frac{\alpha}{\sqrt{n}})\endaligned\tag 2.6
$$
by (2.5).

Next, we expand $L_\rho(\frac{\alpha}{\sqrt{n}})$ around zero, to get
$$
L_\rho(\frac{\alpha}{\sqrt{n}})= 1+\frac{\alpha^2}{2n} +
\frac{1}{6}\frac{\alpha^3}{n^{3/2}}L^{(3)}_\rho (b)
$$
for some $0\leq b\leq \frac{\alpha}{\sqrt{n}}$.
$$
|L^{(3)}_\rho(b) | = | \int x^3 e^{bx}\rho(x) dx |\leq C_\varepsilon
L_\rho(b+\varepsilon)
$$
and we always choose $b+\varepsilon\leq A.$

Finally,
$$
L_{\rho_n}(\alpha) = (1+\frac{\alpha^2}{2n} + \overline{c}_\varepsilon
n^{-3/2})^n \to e^{\alpha^2/2} \text{ as } n\to \infty
$$
so $L_{\rho_n}(a) \leq 2e^{\alpha^2/2}$ for all $n$.\qed
\proclaim{Theorem 2.2 (Polynomial Localization)}  Assume for $N_0$
fixed, $N_0 > 2$
$$
\int |x|^{N_0} \rho(x) dx\equiv M_{N_0}(\rho) < d_0 < \infty.\tag 2.7
$$
Let $\rho_n$ be the normalized $n$-convolution as before.

Then, there exists $d>0$ such that
$$
M_{N_0}(\rho_n)< d(N_0, d_0), \text{ uniformly in } n.\tag 2.8
$$
\endproclaim
\remark{Remark} Similar results with weak localization were proved in
[CS]; they are optimal in the conditions of localization, where
Lindenberg type condition is used.  The proof for such weak
localization is more involved.
\endremark
\demo{Proof}  Consider first $N_0 = 2k, k$ integer.  It is enough to
consider the even case of distribution.

So let $k=2, \quad  \rho,\eta$ even:
$$
\aligned
I_{4,\theta}&=\int x^4\eta(\cos \theta x + \sin \theta y)\rho(-\sin
\theta x + \cos \theta y) dxdy\\
&=\int(u \cos \theta + v \sin \theta)^4 \eta(u)\rho(v) dudv =
\cos^4\theta M_4(\rho) + 6\cos^2\theta \sin^2\theta\\
&+\sin^4\theta  M_4(\rho)\endaligned\tag 2.9
$$
where we used evenness, and the fact that $M_2(\rho)= M_2(\eta) = \int
x^2\rho = \int x^2\eta =1$. (Recall that we always assume that
$M_2(\rho) = 1)$.

Completing to squares, we get from (2.9):

$$
I_{4,\theta} = (M_4^{1/2}(\eta) \cos^2\theta +
M_4^{1/2}(\rho)\sin^2\theta)^2 +
2(3-M_4^{1/2}(\rho)M_4^{1/2}(\eta))\cos^2\theta\sin^2\theta.\tag 2.10
$$

For the gaussian distribution $g$
$$
M_4(g) = 3
$$
Therefore, if $M_4(\rho), M_4(\eta) \leq 3$ the $M_4$ moment increases
under convolution to approach 3.

On the other hand, if both $M_4$ are larger than 3, then
$$
M_4(\rho_2)\text{ decreases, so}
$$
$$
M_4(\rho_2) < \max \{ M_4(\rho), M_4(\eta)\}. \tag 2.11
$$
By Jensen's inequality
$$
M_4(\rho) \geq \left(\int x^2\rho dx\right)^2 \geq 1
$$
so that
$$
M_4^{1/2}(\rho) M_4^{1/2}(\eta) \geq \min\{ M_4^{1/2}(\rho),
M_4^{1/2}(\eta)\}
$$
and hence
$$
3-M_4^{1/2}(\rho) M_4^{1/2}(\eta) \geq 0\text{ only if}
$$
$$
\max\{ M_4^{1/2}(\rho), M_4^{1/2}(\eta)\}\leq 3.
$$
We conclude that
$$
\aligned
M_4(\rho_2)&\leq \max\{ M_4(\rho), M_4(\eta), 9\}\\
\rho_2&\equiv \rho \underset{\theta}\to{*} \eta.
\endaligned
$$
After iteration, we therefore get
$$
M_4(\rho_n) \leq \max\{ M_4(\rho), M_4(\eta), 9\}.
$$
In the case $k> 2$, arbitrary we have in a similar way
$$
\int\int x^{2k}\eta(\theta)\rho(\theta) dxdy =
M_{2k}(\rho)\cos^{2k}\theta + M_{2k}(\eta) \sin^{2k}\theta + R_k
$$
where $R_k$ are lower order moments (in powers of $k$).
And as before, we estimate the above equality by
$$
\leq (C_{2k}\cos^{2}\theta + C_{2k} \sin^2\theta) = C_{2k}
$$
with
$$
C_{2k} \equiv \max\{ M_{2k}(\rho), M_{2k}(\eta), C_{2k-1}\},
$$
from which the result follows.

The general case now follows from the following Proposition (2.3) \qed
\enddemo
\definition{Definition} For a random variable $X$, we define the
$\psi$-function of $X$ as
$$
\psi(R) = E 1_{\{X\geq R\}}X^2 = \int_{|x|\geq R} x^2\rho(x) dx
$$
$E$-expectation, $\bold 1_{\{A\}}$ is indicator function of $A$.
\enddefinition

\proclaim{Proposition 2.3}  Let $\{X_j\}^\infty_{j=1}$ be an
i.i.d. sequence of random variables with $p$ finite
moments, uniformly in $j$, in the integral sense:
$$
\psi_j(R) \leq \psi (R)
$$
and
$$
\int^\infty_1\psi(R)R^{p-3} dR < C_\psi < \infty.
$$
Here $\psi_j(R)$ is the $\psi$-function of $X_j$.

Then, for any $\varepsilon > 0$, there exists a constant $C$,
depending only on $C_\psi$ and $\varepsilon$ such
that
$$
\langle |Z_{2^n}|^{p-\varepsilon}\rangle \leq C(C_\psi,
\varepsilon)\tag 2.12
$$
where
$$
Z_{2^n}\equiv \frac{1}{\sigma_1}\sum^{2^{n-1}}_{j=1} X_j +
\frac{1}{\sigma_2}\sum^{2^{n-1}}_{j=1} X_{j+2^{n-1}}.
$$
\endproclaim

\demo{Proof}  We prove it only for the normalized case where all
variances are 1.

Let $2k< p < 2k + 2$ be given.
$$
Z_{2^n} = 2^{-n/2}\sum^{2^n}_{j=1} X_j =
2^{-n/2}\left(\sum^{2^n}_{j=1} U_j + \sum^{2^n}_{j=1} V_j
\right)\tag 2.13
$$
with
$$
U_j = X_j - V_j
$$
$$
V_j = X_j \bold 1_{\{X_j \leq K\}}.
$$
Then
$$
\aligned \psi_{Z_{2^n}}(R) &= E \bold 1_{\{Z_{2^n}\geq R\}}
Z^2_{2^n} \leq E \bold
1_{\{Z_{2^n}\geq R\}}\left[ 2\{2^{-n/2} \sum U_j\}^2\right.\\
 &+ 2\{2^{-n/2}\left. \sum V_j\}^2\right].\endaligned\tag 2.14
$$
The second term on the r.h.s. of (2.14) is bounded by $2\psi(K)$ and
the first term is controlled by H\"older's inequality:

$$
\aligned \text{ first term } &\leq P(|Z|\geq R)^{\frac{k}{k+1}}
(M_{2k+2}(\tilde U))^{\frac{1}{k+1}} \\
\tilde U &\equiv 2^{-n/2} \sum U_j\\
M_{2k+2} (\tilde U) &\leq C M_{2k+2} (U_1) \leq
\bar{C}K^{2k+2-p}\endaligned
$$
by the even case, where $\bar {C}$ is the p-th moment of $U_1$
$$
P(|Z|\geq R) \leq R^{-2} \psi(R).
$$
Combining all this we get
$$
\psi_{Z_{2^n}}(R) \leq CR^{-2k/(k+1)}\psi (R)^{\frac{k}{k+1}}
K^{2-{\frac{p}{k+1}}} + 2\psi (K)\tag 2.15
$$
Now, choose $K=R$ in (2.15), to get
$$
\psi_{Z_{2^n}} (R) \leq CR^{\frac{2-p}{k+1}} \psi
(R)^{\frac{k}{k+1}} + 2\psi(R).
$$
Multiplying by $R^{p-3-\varepsilon}$ and using H\"older's inequality
again, the result follows.\qed
\enddemo

\subhead Section 3. Propagation of Localization II -
Gaussian\endsubhead

Now we assume that $\rho$ is gaussian localized, normalized
distribution:
$$
\aligned
\int\rho d x &=1 = \int x^2\rho dx\\
|e^{cx^2}\rho(x)|&< C_0 \text{ for some } c> 0, |x|\to
\infty.\endaligned
$$
We use * to denote convolution and $\circledast$ to denote the
normalized (rescaled) convolution: $\circledast = \overset{*}\to{\sqrt{2}}$.

\proclaim{Theorem 3.1}
Let $\rho$ be as above and assume furthermore that
$$
\rho = gF\tag 3.1
$$
and $F$ is logconcave ($ln F$ is concave).

Then $\rho_n=\sqrt{n}\rho *\cdots \rho(\sqrt{n} x)$ is gaussian
localized, uniformly in $n$.
\endproclaim

\demo{Proof}  By Brascamp-Lieb we have that:
$$
\aligned
\rho\circledast \rho &= g F_2\\
\rho\circledast \rho &=\int g(\frac{x+y}{\sqrt{2}}) g
(\frac{x-y}{\sqrt{2}}) F (\frac{x+y}{\sqrt{2}})F(\frac{x-y}{\sqrt{2}})
dy\\
 &= g(x)\int g(y) F(\frac{x+y}{\sqrt{2}})F(\frac{x-y}{\sqrt{2}}) dy = g
F_2\endaligned
\tag 3.2
$$
with $F_2$ logconcave.
\enddemo
Next, we need the following proposition

\proclaim{Proposition 3.2  (Brascamp-Lieb)}

For $g$ Gaussian,
$$
\int x^{2m} gF dx \leq \int x^{2m} gdx \tag 3.3
$$
when $\int gF dx = 1, \text{ and } F$ logconcave.
\endproclaim

>From this proposition it follows that
$$
\int e^{\beta x^2}gF dx \leq \int e^{\beta x^2} gdx\tag 3.4
$$
Since in our case $\rho_n = gF_n$, we get
$$
\rho_n^2 = g^2 F^2_n = (\int g^2 F^2_n dx) (\int g^2 F_n^2 dx)^{-1} g^2
F^2_n = \| \rho_n\|^2_{L^2} g^2 F
$$
$F$ logconcave (since $F_n$ is logconcave).

Hence,
$$
\int e^{\beta x^2}\rho^2_n dx \leq \|\rho_n\|^2_{L^2}\int e^{\beta
x^2} g^2 dx < \infty.
$$
\qed

\remark{Remark}  If $\rho$ is regularized as $\rho\to\rho_t\equiv \rho
\circledast g_t$  we have
$$
\int e^{\beta x^2}\rho_{t,n} = \int \rho_n\circledast g_t e^{\beta
 x^2} = \int \rho_n g_t * e^{\beta x^2} = \int \rho_n e^{\beta_t
x^2}\tag 3.6
$$
with $\beta_t \sim \beta.$

It remains to show that, sufficiently smooth gaussian localized $\rho$, will have the form
$gF$ after sufficiently many iterations.

Next, we demonstrate such cases:
\endremark

\proclaim{Theorem 3.3}  Let
$$\rho=(2\pi)^{-1/2} \exp(-x^2 /2) + p(x)\tag 3.7
$$
and assume that
$$
|\int e^{\alpha x} p(x) dx | \leq C_1 e^{|\alpha|^{2-\varepsilon}},
\varepsilon > 0, \tag 3.8
$$
and $p$ smooth.

Then, for $n$ sufficiently large, $\rho_n = gF$ with $F$ logconcave.
\endproclaim

\demo{Proof} Let, as before
$$\aligned
\rho_\alpha & = \lambda(\alpha)^{-1} e^{\alpha x} \rho(x)\\
\lambda (\alpha) & = \int e^{\alpha x} \rho(x) dx\endaligned
$$
we have a lower bound on $\lambda(\alpha)$:
$$
\lambda(\alpha) = e^{\alpha^2/2} + \int e^{\alpha x} p(x) dx
$$
so, by (3.8) it follows that
$$
\lambda(\alpha) \geq \frac 12(e^{\alpha^2/2} - c) \text{ for }
\alpha > \alpha_0 (C_1, \varepsilon)\tag 3.9
$$
where $\alpha_0$ is approximately $(ln c_1)^\beta$, some $\beta > 0$.

Now,
$$
\aligned
\lambda(\alpha)^{-1}&\int(x-m_\alpha)^4 e^{\alpha x} \rho(x) dx =
\lambda(\alpha)^{-1}\int(x-m_\alpha)^4(2\pi)^{-1/2}
e^{-x^{2/2}}e^{\alpha x} dx\\
&+\lambda(\alpha)^{-1} \int(x-m_\alpha)^4e^{\alpha x} p(x) dx\\
&\equiv I_1 + I_2\endaligned\tag 3.10
$$

$$
\aligned
I_1 &=(2\pi)^{-1/2}\int \{ |x-\alpha|^4 + 6|x-\alpha|^2(m_\alpha -
\alpha)^2 + (m_\alpha-\alpha)^4 + \text{ odd terms
}\}\lambda(\alpha)^{-1} e^{\alpha^{2/2}}\\
&\times e^{-{\frac{1}{2}}(x - \alpha)^2} dx\\
&\leq\{ 3 + 6(m_\alpha - \alpha)^2 + (m_\alpha-\alpha)^4 + 0
\}2e^{\alpha^2/2}/(e^{\alpha^2/2} - C_1).\endaligned\tag 3.11
$$

Therefore $I_1$ remains bounded uniformly in $\alpha$, if $|m_\alpha -
\alpha| \leq C_0$ uniformly in $\alpha$.  Furthermore, $I_2$ is small
when $\alpha$ is large, by our assumptions on $p(x)$.

Now,
$$
\aligned
m_\alpha=\lambda(\alpha)^{-1}\int x e^{\alpha x}\rho(x) dx &= \alpha +
\lambda(\alpha)^{-1}\int x e^{\alpha x} p(x) dx\\
   &=\alpha + O(\alpha^{-\varepsilon})
\endaligned
$$
which implies that the r.h.s of (3.11) is uniformly bounded in
$\alpha$.  To conclude, 3.10 - 3.11 implies that the fourth moment is
uniformly bounded; and the second moment is close to 1.

Next,
$$
\frac{d}{dx}\rho_\alpha = \rho'_\alpha = \text{ nice } +
\lambda(\alpha)^{-1} e^{\alpha x} (\alpha p + p'(x) )
$$
where nice stands for terms which are uniformly bounded in $\alpha$, so,
$$
\|\rho'_\alpha \|^2_{L^2} \leq \| \text{ nice } \|^2 +
\lambda(\alpha)^{-2}\|e^{\alpha x} (\alpha p + p'(x))\|^2_{L^2}.\tag
3.12
$$
$$
\int|e^{\alpha x} (\alpha p + p')|^2 = \int e^{2\alpha x}(p')^2 dx -
\int e^{2\alpha x} \alpha^2 p^2 dx
$$
so, to prove uniformly of a bound on (3.12), in $\alpha$, we only need
to bound
$$
\lambda(\alpha)^{-2}\int e^{2\alpha x}(p')^2 dx \leq C, \text{ uniformly in
} \alpha,
$$
which is implied by our conditions on $p$.

Now, taking the n-th normalized convolution of $\rho_\alpha,
\rho^{(n)}_\alpha$ we know by the polynomial propagation of
localization, Thm 2.1, and by the entropy production bounds of [CS]
that
$$
S(\rho^{(n+1)}_\alpha) - S(\rho^{(n)}_\alpha) > \Phi (S(g_\alpha) -
S(g^{(n)}_\alpha)).
$$
We use that convolution improves or preserves the smoothness of
$\rho$, therefore we can take $\rho$ to be independent of $\rho_u$.
see [CS]:  The function $\Phi_\rho$ was obtained thorough a
compactness argument, and was not computable.  On the other hand, we
were able to show that $\Phi_\rho(t)$ was strictly increasing as a
function at $t$, and hence $\Phi_\rho(t) > 0$ data $\rho$.  Moreover
$\Phi_\rho(t)$ depended on $\rho$ only in a way that was invariant
under the convolution map, so that the same function $\Phi$ could be
used at each stage in the treated convolution.  This act was crucial
in our application which requires us to absorb the effect of
dependence.

In this paper we will estimate $\Phi_\rho$.  We will place more
restrictive conditions on $\rho$, but shall obtain quantitative
information on $\Phi_\rho$ in return.

Hence, $\rho^{(n)}_\alpha$ converges to a gaussian in entropy, $S$,
and so in $L^1$.  By smoothness, all derivatives also converge,
uniformly in $\alpha$.

Now, it follows that for $n> n_0$.
$$
- (ln \rho^{(n)}_\alpha)^{''} |_{x=0} \geq 1-\varepsilon
$$
and since, moreover $\alpha\to m_\alpha$ covers $\Bbb R$, we have that

$$
- (ln \rho^{(n+1)})^{''} \geq 1-\varepsilon\text{ for all } x.
$$

Hence,
$$
\rho^{(n+1)} = \rho_{n+ 1} = e^{-(1-\varepsilon) x^2/2} F
$$
with $F$ logconcave.   \qed
\enddemo

\remark{Remark}  If $\rho$ is not smooth, then we apply the theorems
to $\rho=M\circledast \rho$ with $M$ gaussian.  For such $\rho$ the
condition on $p'$ is satisfied whenever we have the bound 3.8, since
$$
p'=M'\circledast p.
$$

Furthermore, the gaussian localization of $(M\circledast \rho)_n$
implies that of $\rho_n$, since
$$
(M\circledast\rho)\circledast (M \circledast \rho) = M
\circledast(\rho\circledast\rho)
$$
so, since $M$ is well localized, $\rho\circledast\rho\circledast \cdots \rho$
is well localized whenever
$$
(M\circledast\rho)\otimes(M\circledast \rho) \cdots (M \circledast \rho)\text{ is
well localized.}
$$
\endremark

\subhead Section 4. Entropy Production
\endsubhead

In this section, we prove optimal entropy production bounds for the
convolution map.

Recall the following formula for the Entropy production by convolution
[CS]
$$
S(\rho\circledast\rho) - S(\rho) = \int^\infty_0 J(\rho_t\circledast
\rho_t) - J(\rho_t) dt\tag 4.1
$$
where $S$ is the entropy and $J$ is the relative information.

$\rho_t$ is the map, up to time $t$ of $\rho$ under the Orenstein-Uhlenbek
process.

Also from [CS, Bar] we have the following bounds
$$
|\nabla\sqrt{\rho(x)} |^2 \leq B_t P^*_t \rho(x) \tag 4.2
$$
which, by the way of the localization of $\rho$ implies that
$|\nabla\sqrt{\rho}|^2$ is similarly localized.

Also, recall the definition of the $\psi$ function
$$
\psi(R)_ =\int_{|x|\geq R} x^2\bold \rho(x) dx.
$$

Define
$$
J_R(\rho) = 4\int_{|x|\geq R} |(\nabla + \frac{x}{2})
\sqrt{\rho}|^2 dx.
$$

\proclaim{Lemma 4.1}
$$J_R(\rho) \leq 2\psi(R) + 8(1+ R^2)^{-1} B_t P^*_t\psi(R).\tag
4.3
$$
\endproclaim

\demo{Proof}  Follows from (4.2) and the definition of $\psi(R)$.
\enddemo

\proclaim{Lemma 4.2}
$$P^*_t\psi(R)\leq \psi_\rho(R/2)+\psi_g(R/2)\tag 4.4
$$
\endproclaim

\demo{Proof}
See [CS]
\enddemo

We can now state the main entropy production bound : (see CC1, CS for similar results
with weaker nonlinear (lower bounds) in $D(\rho)$, in the case of
Boltzman equation and the CLT, respectively.  However, those results
do hold for general $\rho$; i.e. finite variance and finite entropy
are the only conditions imposed.)

\proclaim{Theorem 4.3} Let $\rho$ satisfy $J(\rho), S(\rho)$
finite, $\rho$ smooth, and have a finite second moment.

\noindent{\it (1)}  Suppose that $K\ge g/\rho \ge 1/K$ for some constant $K$. Then
$$
S(\rho\circledast \rho) \geq \frac{K}{2} D(\rho).\tag 4.6
$$

\noindent{\it (2)}  More generally, define $R_\epsilon$ so that
$$ 2\psi_\rho(R_\epsilon) + 8(1+ R_\epsilon ^2)^{-1}  + \psi_\rho(R_\epsilon/2)+\psi_g(R_\epsilon/2)
< J(\rho)/2 := \epsilon\ .$$

 Suppose that $g/\rho$ is bounded below by by $K_\epsilon$ on the ball of radius $R_\epsilon$. Then
$$
S(\rho\circledast \rho) \geq C_\varepsilon D(\rho).\tag 4.6b
$$
where $C_\epsilon$ depends only on $\epsilon$ and $\psi_\rho$.
\endproclaim

\proclaim{Remark} The constant $C_\varepsilon$ depends on the
localization of the relative Fisher information, and the distance of
the distribution $\rho$ from the normalized Gaussian. Therefore, an
estimate with known, uniformly bounded constant, would require
controlling such quantities. This follows when we have propagation
of Gaussian localization, as in Section 3. Alternatively, one may
expect to prove propagation of localization for the relative Fisher
information, which we do not have. In Section 5, we use a new
construction (stitching), to obtain uniform bounds for
$C_\varepsilon$.
\endproclaim

 \demo{Proof}

If $\rho= g$ there is nothing to prove.

For $\rho\neq g, J(\rho) > 0$.  So assume $J(\rho) =
\varepsilon$.  We now choose $R$ so large that
$$J_R(\rho) \leq \frac{1}{2} J(\rho)
$$
$R(\varepsilon)$ is fixed by
$$
2\psi(R) + C_t\tilde\psi(R)/(1+R^2) \leq \varepsilon/2 \tag 4.7
$$
with $\tilde\psi \equiv P_t^* \psi$.

Next, we use the lower bound, proposition (4.4) below:
\medskip

$$
\aligned
J(\rho_t) &- J(\rho_t\circledast\rho_t)\geq F_{a,\tau}\\
&\equiv\inf_{c,d}\{ E\left[\frac{d}{dx} \ln \rho_t(\frac{\tau}{a} G)+cG
+ d\right]^2\}\\
&=\inf_{c,d}\int | \nabla \ln \rho_t(x) + c x + d|^2 g(x) dx\\
&\geq \int_{|x|\leq R(\varepsilon)} |\nabla \ln \rho_t(x) -c^*
x-d^*|^2 g(x) dx\\
&\text{for some } c*, d^*.\\
&\text{This last expression is then equal to}\\
&=\int_{|x|\leq R(\varepsilon)} |Q |^2 \frac{g(x)}{\rho_t(x)} \rho
_t (x) dx \\
&\geq \int_{|x|\leq R(\varepsilon)} |Q |^2\rho_t(x) dx
\cdot \| \frac{\rho_t (x)}{g(x)} \|^{-1}_{L^\infty(|x|\leq
R(\varepsilon))}\\
&\text{with } Q=\nabla\ln \rho_t(x) - C^*x-d^*,
\endaligned\tag 4.8
$$
and we also have
$$
\int_{|x|> R(\varepsilon)} |\nabla\ln \rho_t(x) - c^* x - d^*|^2\rho(x)dx
\leq \varepsilon/2.\tag 4.9
$$

Finally, (4.8) and (4.9) imply
$$
\aligned
J(\rho_t) - J(\rho_t\circledast \rho_t) &\geq \|
 \, \frac{\rho}{g}\|^{-1}_{L^\infty(|x|\leq R(\varepsilon))}\frac12
\int|\nabla ln \rho_t-c^*x-d^*|^2\rho_tdx\\
&\geq \frac{1}{2}\|\frac{\rho}{g}\|^{-1}_{L^\infty(R(\varepsilon))} J
(\rho_t).\endaligned\tag 4.10
$$
The theorem now follows from this last inequality and (4.1).\qed
\enddemo

\proclaim{Proposition 4.4}
$$
J(\rho_t) - J(\rho_t * \rho_t) \geq F_{a,\tau} \equiv
\inf_{c,d}\{\int|\frac{d}{dx} ln \rho_t(x) + c x+ d|^2 g(x) dx\}.
\tag 4.11
$$
\endproclaim
\demo{Proof}  Introduce the convolution operator $C_{\rho,\theta}$
$$
C_{\rho,\theta} f \equiv \int f(\langle e_1, R_\theta(x, y)\rangle)
\rho(y) dy
$$
where $\langle , \rangle$ is the scalar product in $\Bbb R^2, e_1=(1,0)$
and $R_\theta$ is rotation in $\Bbb R^2$ by $\theta$.
$$
C_{\rho,\theta}:L^2(\rho) \to L^2(\rho)
$$
for any $\rho= g, g$ gaussian and $\theta = e^{-t},\quad C_{\rho,
t}$ becomes the Orenstein-Uhlenbek process.

In this case $C_{\rho,\theta}$ is self-adjoint and its eigenvalues are
$\cos^n\theta$.

In general $C_{\rho,\theta}$ is not bound on $L^2$ and is selfadjoint
only for $\rho = g$.

Let $\Pi_j$ denote the projection on the subspace of the first $j$
eigenvectors of $C_{\rho,\theta}$.
$$
\Pi_j + \bar{\Pi}_j = \bold 1.
$$
Now, consider
$$I_\theta \equiv \int\int | h(x) + h(y) -
\bar{h}(R_\theta(x,y))|^2\rho(x) \rho(y) dx dy.
$$
The following lemma is essentially due to Brown [Br].  See [CC2] for
an adaptation to the Boltzmann equation setting.
\enddemo
\proclaim{Lemma 4.5 (Linear Approximation Lemma)}
$$
I_\theta\geq C_\theta\inf_{a,b} \int |h(x) -ax-b|^2 \rho(x) dx.\tag
4.12
$$
\endproclaim

See [Br]. Here we use it with $ \theta=\pi/4.$

 \subhead Section 5. How to deal with thin tails \endsubhead
\medskip

\proclaim{Lemma 5.1} Let $\rho$ be a probability density with $I(\rho) < \infty$.  Then
for $q>1$ and $R>0$,
$$\int_{\{|x|> R\}}\rho^q(x)dx \leq  I(\rho)^{q-1}\left(\int_{\{|x|> R\}}\rho(x)dx\right) \ .$$
\endproclaim

\noindent{\bf Proof:} Let $f := \sqrt{\rho}$. Using the bound $\|f\|_\infty^2 \leq 2\|f\|_2\|\nabla f\|_2$ for functions on $R$,
$$\int_{\{|x|> R\}}\rho^q(x)dx = \int_{\{|x|> R\}}f^2 f^{2(q-1)}(x)dx \leq
\left(\int_{\{|x|> R\}}\rho(x)dx\right)(2\|\nabla f\|_2)^{2(q-1)}\ .$$
Recall  that $2 \|\nabla f\|_2 = \sqrt{I(\rho)}$. \qed
\medskip

\proclaim{Lemma 5.2} Let $\rho$ be a probability density with $I(\rho) < \infty$ and finite second moment.
Then
$$
\int_{|x|\geq R}\rho |\ln \rho| dx \leq
2I(\rho)^{1/2}  \left(\int_{\{|x|> R\}}\rho(x)dx\right) + \frac{ \sqrt{\pi}}{2}\left(\int_{|x|>R} \rho (1+|x|^2)dx\right)^{1/2}\ .
$$
\endproclaim

 \noindent{\bf Proof:}
Fix any $r>0$. On the set $\{ \rho > 1\}$,
$$\rho |\ln \rho| = \rho \ln {\rho} \leq \frac{1}{r}(\rho^{1+r} - \rho) \leq \frac{1}{r}\rho^{r+1}\ .$$
By the previous lemma,
$$\int_{\{\rho \geq 1\}\cap\{|x|\geq R\}} \rho |\ln \rho|  \leq \frac{1}{r}I(\rho)^r  \left(\int_{\{|x|> R\}}\rho(x)dx\right) \ .$$

On the set $\{ \rho < 1\}$,
$$\rho |\ln \rho| = \rho \ln \frac{1}{\rho} \leq \frac{1}{r}(\rho^{1-r} - \rho) \leq \frac{1}{r}\rho^{1-r}\ .$$
Therefore, by H\"older,
$$\eqalign{
\int_{\{\rho \leq 1\}\cap\{|x|\geq R\}} \rho |\ln \rho|  &\leq \frac{1}{r} \int_{|x|\geq R}  \rho^{1-r} <x><x>^{-1} dx\cr
& \leq
\frac{1}{r}\left(\int_{|x|\geq R}\rho<x>^{1/(1-r)}d x\right)^{1-r} \left(\int <x>^{-1/r} \right)^r\ .\cr}$$
Choosing $r=1/2$, we obtain the result. \qed

\proclaim{Proposition 5.3}  Let $\rho$ be a probability density mean zero, unit variance, $I(\rho) < \infty$ and finite third moment.
Let
$$\rho_n= \rho\circledast\rho\cdots \circledast \rho  \quad  n -  \text{times}\ .$$
Then there exists a constant $c$ such that for all $n$,
$$\int_{|x| < R}|\rho_n - g|d x \leq c R 2^{-n/2}\ .$$
$$\int_{|x| < R}|\rho_n/g - 1|d x \leq c Re^{R^2/2} 2^{-n/2}\ .$$
\endproclaim

\noindent{\bf Proof:} See Feller or Major

We are now ready to define the stitching operations.

Recall the definition
$$
\rho_{2n}:= \sqrt2\rho_{n-1}*\rho_{n-1}(\sqrt2 x),
$$
with $\int x^2\rho_n=\int x^2\rho_0=1.$ We further define $N:=2^n.$
Then, we let, for some fixed $c>0,$
$$
\tilde{\tilde\rho_n}:=\rho_n \chi_{c\sqrt
n}+\frac{1}{\sqrt{2\pi}}e^{-\frac{x^2}{2}} (1-\chi_{c\sqrt n}),
$$
where
$$
\chi_m:= h_0*I_{[-m,m]},
$$
with a nonnegative mollifier function $h_0$, satisfying: $h_0\geq
0$, $h_0\in C_0^{\infty}$, Support of $h_0\in [-1,1]$, $\int h_0=1.$
Here $I_B$ denotes the characteristic function of the set $B.$ We
then normalize :
$$
\tilde \rho_n(x):=c_n\tilde {\tilde\rho_n}(d_nx-e_n),
$$
such that $$ \int \tilde \rho_n=1, \int x^2\tilde \rho_n=1,$$

$$\int x \tilde \rho_n=0.$$

 Writing $c_n=1+\epsilon_n, d_n=1+\epsilon_n^{\prime}, e_n=1+\epsilon_n^{\prime \prime}$, it
follows, by an application of the local central Limit Theorem, and
localization, that the $\epsilon_n$'s tend to zero, as $n$ goes to
infinity.

 \proclaim{Proposition 5.4}
 Let $S$ denote the entropy functional, as before, and $\rho_n,
 \tilde\rho_n$ defined as above.
 Then,
 $$
 S(\rho_n)-S(\tilde \rho_n)= r(n^{1/2})^{-1}.
 $$
$r(k)$ tends to infinity as $k$ goes to infinity. Moreover, if
$\rho_1$ is polynomially localized to order $2m+2$, then $r(k)$
grows like $k^m$; for $\rho_1$ exponentially localized, $r(k)$ is
exponentially growing in $k$.
\endproclaim

\noindent{\bf Proof:}
$$
S(\rho_n)-S(\tilde \rho_n)=\int_{|x|\leq c\sqrt
n}(\rho_n\ln\rho_n-\tilde \rho_n \ln \tilde \rho_n)+ R_n
$$
where $$ R_n=\int_{|x|\geq c\sqrt n}(\rho_n\ln\rho_n-\tilde \rho_n
\ln\tilde \rho_n).
$$
If $\rho_1$ is polynomially localized, to order $2m$, (respectively,
exponentially localized), then by our previous results on
propagation of localization, in these cases, the localization
persists, uniformly in $n. $ Since the range of integration in the
$R_n$ term is $|x|\geq c\sqrt n$, the bound $R_n=r(n^{1/2})^{-1}$
follows.

It remains to control the other part of the integration region. In
this region we have that:
$$
\tilde \rho_n= \frac{\rho_n}{1+\epsilon_n},
$$
and therefore,
$$
\aligned
 \int_{|x|\leq c\sqrt n}(\rho_n\ln\rho_n-\tilde \rho_n \ln
\tilde \rho_n)\\
&= \int_{|x|\leq c\sqrt
n}(\rho_n\ln\rho_n-(1+\epsilon_n)^{-1}\rho_n)[\ln \rho_n-\ln(
1+\epsilon_n)]\\
&= \int_{|x|\leq c\sqrt n}\rho_n\ln\rho_n
(1-\frac{1}{1+\epsilon_n})\\&+\int_{|x|\leq c\sqrt
n}(1+\epsilon_n)^{-1}\rho_n\ln(1+\epsilon_n).
\endaligned
$$

Since the entropy is uniformly bounded in $n$, and the $\rho_n$ are
all normalized to 1, the proof follows, if we show that $$
\epsilon_n=r(n^{1/2})^{-1}.$$

This last estimate follows directly from the definition of the
stitched distribution:
$$
\int \tilde {\tilde\rho_n}=\int_{|x|\leq c\sqrt n}\rho_n+
R_n=R_n+R_n+1.
$$
 Similar estimate holds for for the other $\epsilon$'s.

 \qed
\proclaim{Proposition 5.5}

Let $N$ be defined as before, for any fixed $n.$ Assume that $\rho$
satisfies the normalization conditions as before, and furthermore it
is Gaussian, exponential or polynomially (of order $p\ge4$) localized:
$$
\|e^{bx^2}\rho\|_{\infty} \leq 1, \quad \quad b>0.
$$

$$
\|e^{b|x|}\rho\|_{\infty} \leq 1, \quad \quad b>0.
$$

$$
\||x|^p\rho\|_1 \leq b, \quad \quad b>0.
$$

Let $\tilde \rho$ be the associated stitched distribution as defined
before. Then,
$$
S(\rho_{2N})\ge S(\tilde\rho_N*\tilde\rho_N)-r(\sqrt N)^{-1}.
$$
\endproclaim
\noindent{\bf Proof:}

Using that
$$
S(\rho)= \sup_{\phi\in O} \left( \int\rho \phi dx- \ln\int e^{\phi}
dx\right ),
$$
and choosing $ e^{\phi}=\tilde\rho_n*\tilde\rho_n,$ n arbitrary, we
arrive at:
$$
\align
 S(\rho_{2n})>S(\tilde\rho_n*\tilde\rho_n)+ \int \rho_{2n}\ln (\tilde\rho_n*\tilde\rho_n)\\
 &=\int_{|x|<c\sqrt n/2} (\rho_{2n}-\tilde\rho_n*\tilde\rho_n)\ln
(\tilde\rho_n*\tilde\rho_n)\\
&+\int_{|x|>2c\sqrt n/2} (\rho_{2n}-\tilde\rho_n*\tilde\rho_n)\ln
(\tilde\rho_n*\tilde\rho_n)\\
&+\int_{c\sqrt n/2<|x|<2c\sqrt n}
(\rho_{2n}-\tilde\rho_n*\tilde\rho_n)\ln
(\tilde\rho_n*\tilde\rho_n)\\
&= S(\tilde\rho_n*\tilde\rho_n)+0-\int_{|x|>2c\sqrt n}
(x^2/2)(\rho_{2n}-\tilde\rho_n*\tilde\rho_n)+ B,
\endalign
$$
$$
B:=  \int_{c\sqrt n/2<|x|<2c\sqrt n}
(\rho_{2n}-\tilde\rho_n*\tilde\rho_n)\ln
(\tilde\rho_n*\tilde\rho_n).
$$
We now use this last inequality with $n$ replaced by $N:=2^n.$

 Then,
we choose $2c<c_0$, so that for $n>N_0$, we have that $\rho_N\ge
e^{-x^2/3}$ for $ |x|\leq 2c\sqrt N.$ Hence
$$
B \leq CN\int_{2c>|x|>c\sqrt N/2}
(\rho_{2N}-\tilde\rho_N*\tilde\rho_N)\leq c_1 e^{-cN}
$$

since, by the pointwise CLT, for such $x$, we have gaussian localization.

Finally,
$$
-\int_{|x|>2c\sqrt N}
(x^2/2)(\rho_{2N}-\tilde\rho_N*\tilde\rho_N)=r(\sqrt N)^{-1}.
$$\qed

\medskip
\noindent{\bf Proof of the Main Theorem-I}

By the above proposition we have that:
$$
\align
S(\rho_{2N})\geq S(\tilde\rho_N*\tilde\rho_N)-\\r(\sqrt N)^{-1}\\
&\geq S(\tilde\rho_{N})+ \Phi(S(\tilde\rho_N|g))-\\r(\sqrt N)^{-1}
&\geq S(\rho_N)+ \Phi(S(\tilde\rho_N|g))-r(\sqrt N)^{-1}
\endalign
$$
The proof of the main theorem ,namely that $S(\rho_N)\longrightarrow
S(g)+r(\sqrt N)^{-1}$, follows from the following:

 \proclaim{Theorem 5.6}
For $\rho$ Gaussian localized as above, and for all $n$ large
enough, we have:

$$
 c_1g\leq \tilde \rho_n \leq c_2g,
 $$
 $$
0\leq S(\tilde\rho_n*\tilde\rho_n |g)\leq (1-c)S(\tilde\rho_n|g).
$$

$c$ depends on $c_1,c_2$, and $0<c<1.$
\endproclaim

The proof of the above theorem follows from the construction of
$\tilde\rho$ and our previous estimates on entropy production in the
Gaussian localized case.

\noindent{\bf Completion of the Proof of the Main Theorem}

The proof now follows, since we can replace
$\Phi(S(\tilde\rho_N|g))$ by $c(S(\tilde\rho_N|g)),$ $c$ is strictly
positive, uniformly in $N,$ since $c_1,c_2$ can be chosen uniformly
in $N$, for all $N$ large enough.\qed

 Then, the relative entropy satisfies, under favorable localization
conditions
$$
D(\rho_{2N}) - D(\rho_N) \geq \delta_0 D(\rho_N).\tag 5.1
$$
>From this, we immediately conclude that the relative entropy
converges to zero, exponentially fast in $N$.

This is the basis for the argument giving an optimal convergence rate
in the Entropy sense, for localized initial distributions $\rho$.

The  inequality (5.1) is the crucial inequality, proved in sections
3, using the propagation of localization for gaussian localized
$\rho$. The {\bf MAIN THEOREM} now follows:

\demo{Proof}  Since $\rho$ is gaussian (or exponentially or polynomially) localized and smooth, we see
that $\rho$ satisfies the conditions
 for Theorems 5.4,5.5,5.6.

Hence, either (in the gaussian or exponential case)
$$
\int\tilde\rho_n e^{\beta |x|} dx < c<\infty, \text{ independently
of $n$},\tag 5.3a
$$

or,
$$\int\tilde\rho_n  |x|^p dx < c<\infty, \text{ independently
of $n$},\tag 5.3b
$$

Next, we apply Theorems 4.3,5.4-5.6 to $\rho_N$ to conclude that
$$
D(\rho_{2N}) - D (\rho_N) \geq C_\varepsilon D(\rho_N)-r(\sqrt N)^{-1},
\tag 5.4
$$
with
$$
C_\varepsilon = C\|\frac{\tilde\rho_N}{g}
\|^{-1}_{L^\infty(R(\varepsilon))}.\tag 5.5
$$

Due to the propagation of localization (5.3), we see that
$\|\frac{\tilde\rho_N}{g}\|_{L^\infty(R(\varepsilon))}<\infty$,
uniformly in $N$ and hence $C_\varepsilon > \delta>0$ uniformly in
$N$, which implies that
$$
|D(\rho_N)|\leq r(\sqrt N)^{-1} +O(1/N).
$$
\qed\enddemo

\enddocument

\end